%% file: ccc5.tex
\documentstyle[12pt,psfig]{amsart}
\input{def.tex}
\begin{document}
\title{Pseudo-holomorphic curves and the Weinstein conjecture}
\author{Weimin Chen}
\date{\today}
\maketitle
\sectioni{Introduction}
In this paper we will illustrate how certain non-vanishing theorems of Gromov-Witten invariants can be used to prove the Weinstein conjecture.

\vspace{2mm}

Let $S\subset (N,\omega)$ be a hypersurface in a symplectic manifold. The 
characteristic distribution $\E_S$ on $S$ consists of the tangent vectors 
$v\in TS$ such that $i(v)\omega_S=0$, where $\omega_S$ is the pull back of $\omega$ to $S$. The flow lines generated by a vector field in $\E_S$ are 
called characteristics. $S\subset (N,\omega)$ is said to be 
of contact type if there is a 1-form $\alpha$ on $S$ such that $d\alpha=\omega_S$ and $\alpha(v)\neq 0$ for any $v\neq 0$ in $\E_S$.
The 1-form $\alpha$ is called a contact form on $S$, which induces a contact
structure $\xi=\{\alpha=0\}$ on $S$ (see [W]).

\begin{conj}
(Weinstein [W])

If $S$ in $(N,\omega)$ is a compact hypersurface of contact type with
$H^1(S,\R)=0$, then $S$ has at least one closed characteristic.
\end{conj}

We will focus on dimension $4$. Our main result is 

\begin{thm}
Let $S\subset (X,\omega)$ be a compact hypersurface of contact type in a closed symplectic 4-manifold with $b^+_2(X)>1$. Let $\alpha$ be an induced contact
form and $\xi=\{\alpha=0\}$ be the corresponding contact structure on $S$. If 
we orient $S$ by $\alpha\wedge d\alpha$, then 
\begin{enumerate}
\item $c_1(\xi)$ is Poincare dual to a finite union of closed characteristics in $S$,
each of which is oriented by $-\alpha$. In particular, $S$ has at least one closed characteristic if $c_1(\xi)\neq 0$.
\item If $S$ bounds a submanifold $W\subset X$ such that $c_1(W)\neq 0$ and $\omega$ is exact on $W$, then $S$ has at least one closed characteristic.
\end{enumerate}
\end{thm}

As a corollary of Theorem 1.2, we prove the Weinstein
conjecture for compact hypersurfaces of contact type in a 2-dimensional Stein manifold under a mild restriction.

Let $(V,J)$ be a Stein surface with a strictly plurisubharmonic function $\phi$.
Then $\omega_\phi=-d(J^\ast(d\phi))$ is a Kahler form on $(V,J)$ and $(V,\omega_\phi)$ is a symplectic manifold ([EG]). Suppose that $S$ is a
connected orientable compact hypersurface in $V$. Then $S$ bounds a compact submanifold $W\subset V$ since $H_3(V)=0$. 

\begin{coro}
If $S\subset (V,\omega_\phi)$ is of contact type and $c_1(W)\neq 0$, then $S$
has at least one closed characteristic. Moreover, suppose $\alpha$ is an induced
contact form on $S$ and $S$ is oriented by $\alpha\wedge d\alpha$, then the first Chern class of the induced contact structure on $S$ is Poincare dual to a 
finite union of closed characteristics in $S$, each of which is oriented by $-\alpha$. 
\end{coro}

\noindent{\bf Proof:}
Pick a regular value $r>\max(\phi|_{S})$ of $\phi$ so that 
$S\subset\overline{V_r}=\{\phi\leq r\}$ which is compact. Note that $H_3(\overline{V_r})=0$ because $\overline{V_r}$ is homotopic to a 2-complex.
It then follows that
$\overline{V_r}\setminus S$ is disconnected, since otherwise $S$ would have 
intersected a 1-dimensional cycle geometrically once which contradicts the fact that the fundamental class of $S$ (note that $S$ is orientable) in $\overline{V_r}$ is zero. Now $S$ 
must bound a compact submanifold $W\subset V_r$ since $\partial\overline{V_r}$ is connected.

On the other hand, by Theorem 3.2 in [LM], $V_r=\{\phi<r\}$ admits a holomorphic embedding into
a Kahler surface $X$ with $b^+_2(X) >1$ such that the pull back of the Kahler form equals to $\omega_\phi$. It is easy to see that  Corollary 1.3 follows from Theorem 1.2.
\hfill $\Box$

\vspace{3mm}

The Weinstein conjecture has a version for contact manifolds. Let $(M,\alpha)$
be a closed contact manifold with contact form $\alpha$. The Reeb vector field $v$ on $M$ associated to $\alpha$ is defined by $i(v)d\alpha=0$ and $\alpha(v)=1$. One can regard $M=\{0\}\times M\subset\R\times M$ as a compact hypersurface of contact type in the symplectization $(\R\times M,d(e^t\alpha))$ with an induced contact form $\alpha$. The induced characteristics in $M$ are the orbits of the Reeb vector field on $M$.
The Weinstein conjecture for a contact manifold $(M,\alpha)$ states that
$M$ has at least one closed Reeb orbit if $H^1(M,\R)=0$.
 
Let $(Y,\alpha)$ be a closed contact 3-manifold.  Hofer proved in [H] that $Y$ has at least one closed Reeb orbit if $Y=S^3$ or $Y$ is covered by $S^3$; or $\pi_2(Y)\neq 0$; or $\alpha$ induces an overtwisted contact structure on $S$. The remaining case is that $\pi_2(Y)=0$ and $\alpha$ is tight. 

Basic examples of $(W,S)$ in Corollary 1.3 are the Stein surfaces with boundary
([Go]), which have tight induced contact structures on the boundary. So our result is somehow complement to Hofer's. Moreover, the closed Reeb orbits
obtained by Hofer in [H] are all contractible while we can produce closed Reeb orbits of non-trivial homology class when the first Chern class of the contact 
structure is non-zero. However, whether any connected tight contact 3-manifold $(Y,\alpha)$ is fillable, especially Stein fillable, is still a basic open question in contact geometry, although one can construct many concrete examples of that sort ([Go]).

\vspace{3mm}

In the proof of Theorem 1.2 we make use of the following deep theorem of Taubes.

\begin{thm} (Taubes [T])

Let $(X,\omega)$ be a closed symplectic 4-manifold with $b^+_2(X) >1$ and a non-trivial canonical bundle $K$. Then for a generic $\omega$-compatible 
almost complex structure $J$, the Poincare dual to $c_1(K)$ is represented by the fundamental class of an embedded $J$-holomorphic curve $\Sigma$ (may have several components) in $X$.
\end{thm}

The proof starts with the observation that the induced contact structure $\xi$
as a complex line bundle
on $S$ is isomorphic to $K^{-1}|_{S}$ where $K$ is the canonical bundle of 
$(X,\omega)$.
By Theorem 1.4, $c_1(K)$ is Poincare dual to an embedded $J$-holomorphic curve
$\Sigma$ in $X$ for a generic almost complex structure $J$. It is easy to see
that $c_1(\xi)$ is Poincare dual to $\Sigma\bigcap S$ in $S$. The upshot is that
$\Sigma\bigcap S$ converges to a union of closed characteristics in $S$ through a deformation of almost complex structures in a neighborhood of $S$.
The technical part of the proof is a combinatorial argument which produces convergent sequences of annuli in $\R\times S$
from the pseudo-holomorphic curves $\Sigma$ such that the $S$-component of each limiting annulus lies in a Reeb orbit (see Lemmas 2.8, 2.9). The fact that each component of the pseudo-holomorphic curves has a uniformly bounded genus plays an essential
role in the proof. The assumption $c_1(W)\neq 0$ is used to ensure that each pseudo-holomorphic curve goes through $W$, and the assumption that $\omega$ is exact on $W$ is needed to get a non-trivial closed Reeb orbit.

\vspace{3mm}

We end this introduction with two remarks. First, analogous versions of Theorem 1.2 and Corollary 1.3 hold for $b^+_2=1$ case or higher dimensions provided that an analogous non-vanishing theorem of Gromov-Witten invariants exists.
In fact, after this paper was finished, we learned that some relevant results (for any dimensions), especially the stabilized version of the Weinstein conjecture, had been proved recently by G. Liu and G. Tian
in [LT] by a different argument. 
Second, the first Chern class of the contact structure being Poincare dual to a union of closed Reeb orbits is the analogue of Taubes' theorem (Theorem 1.4)
for contact 3-manifolds which can be realized as a hypersurface of contact type
in a closed symplectic 4-manifold with $b^+_2>1$; it would be interesting to know whether this is true for more general contact 3-manifolds, as well as what contact 3-manifolds can be realized as a hypersurface of contact type in a closed symplectic 4-manifold with $b^+_2>1$.

\vspace{3mm}

\centerline{\bf Acknowledgments}

\vspace{2mm}
I am indebted to Yasha Eliashberg and Cliff Taubes for useful communications. I am also grateful to Selman Akbulut for encouragement and many good suggestions after reading the draft version of this paper, to John McCarthy, Tom Parker and Jon Wolfson for their questioning and interest in the work, and to Ron Fintushel who generously supported me (in part) to attend the
1997 IAS/Park City Mathematics Institute Graduate Summer School, during which
this work was initiated.

\sectioni{The proof}
We first recall some basic facts about pseudo-holomorphic curves (or maps)
in (or into) a hermitian manifold. Let $(N,J)$ be an almost complex manifold. A $C^1$ map $f$ from a
Riemann surface $(\Sigma,j)$ into $(N,J)$ is said to be pseudo-holomorphic if the equation $J\circ df=df\circ j$ holds. The image of
$f$ is called a pseudo-holomorphic curve in $(N,J)$. If we choose a Kahler
metric $\mu$ on $(\Sigma,j)$ and a hermitian metric
$h$ on $(N,J)$, the energy of a map $f:\Sigma\rightarrow N$ is 
$$
E(f)=\frac{1}{2}\int_{\Sigma}|df|^2.
$$
Note that the energy $E(f)$ depends only on the complex structure $j$ of $\Sigma$ and the metric $h$ on $N$, while the integrand $|df|^2$ also depends on the choice of the metric $\mu$. For a pseudo-holomorphic map $f:(\Sigma,j)\rightarrow (N,J,h)$, the area of the image of $f$ measured with the metric $h$ is equal to the energy of $f$:
$$
Area_h(f)=E(f)=\frac{1}{2}\int_{\Sigma}|df|^2.
$$

Let $(N,\omega)$ be a symplectic manifold of dimension $2n$. Since $U(n)$ is
a deformation retract of $Sp(2n)$ --- the group of $2n$-dimensional symplectic
matrices,
the tangent bundle of $(N,\omega)$ admits an almost complex structure and any
two such structures are homotopic. An almost complex structure $J$ on $(N,\omega)$ is said to be $\omega$-compatible if $\omega(\cdot,J\cdot)$ is a
hermitian metric on $(N,J)$. It is a well-known fact that the space of all
$\omega$-compatible almost complex structures is nonempty and contractible. 
In a symplectic manifold $(N,\omega)$, for any $\omega$-compatible
almost complex structure $J$ and the associated hermitian metric 
$h=\omega(\cdot,J\cdot)$, we have
$$
Area_h(f)=E(f)=\frac{1}{2}\int_{\Sigma}|df|^2=\int_{\Sigma}f^\ast\omega
$$
for a pseudo-holomorphic map $f:\Sigma\rightarrow N$.
In particular, the pseudo-holomorphic maps (curves) representing a fixed 
homology class in $N$ have uniformly bounded energy (area). In fact, they are
absolutely area minimizing in the homology class, and the images are minimal
surfaces. Any $C^1$ pseudo-holomorphic map into an almost complex manifold
$(N,J)$ is smooth, and a bounded sequence of pseudo-holomorphic maps $f_n:(\Sigma,j,\mu)\rightarrow (N,J,h)$ with a uniform $C^0$ upper bound on the 
gradient $df_n$ (in fact only a uniform $L^p$ bound for $p>2$ is needed) has a 
$C^\infty$ convergent subsequence on any compact subset of $\Sigma$.

We particularly need the following well-known local properties of pseudo-holomorphic maps (curves) into (in) a hermitian manifold (see, for
example, [PW] or [Ye]).

\begin{lem}
(Energy estimate)

Let $(\Sigma,j,\mu)$ be a compact Riemann surface and $(N,J,h)$ be a hermitian manifold with bounded geometry. Suppose $U\subset \Sigma$ is an open subset such that $U\bigcap\partial\Sigma=\emptyset$. Then there exists a constant $\epsilon >0$ such that for any pseudo-holomorphic map $f:\Sigma\rightarrow N$ with
$$
\int_{B(z,2r)}|df|^2 < \epsilon
$$
on a geodesic disc $B(z,2r)$ in $U$ centered at $z\in U$ with radius $2r$, the following estimate holds for a constant $C>0$:
$$
\sup_{B(z,r)}|df|^2\leq\frac{C}{r^2}\int_{B(z,2r)}|df|^2.
$$
\end{lem}

\begin{lem}
(Monotonicity)

Let $(N,J,h)$ be a hermitian manifold with bounded geometry. Then there exist constants $C_0,r_0>0$ with the following property: let $(\Sigma,j)$ be a 
compact Riemann surface with boundary, for any pseudo-holomorphic map $f:(\Sigma,j)\rightarrow (N,J)$, if $f(\partial\Sigma)$ lies outside of a closed $r$-ball $B(f(p),r)$ in $N$ for some $p\in\Sigma\setminus\partial\Sigma$ and $r\leq r_0$, the following inequality holds:
$$
Area_h(f(\Sigma)\bigcap B(f(p),r))\geq C_0 r^2.
$$
\end{lem}

The following lemma is a well-known fact of which we give a short proof here
for completeness.

\begin{lem}
$S\subset (X,\omega)$ is of contact type with a contact form $\alpha$ if and only if a neighborhood of $S$ in $X$ is symplectomorphic to $(-\delta,\delta)\times S$ with the symplectic form
$d(e^t \alpha)$.
\end{lem}

\noindent{\bf Proof:}
The ``if'' part is obvious; we give a proof for the ``only if'' part.

According to [W], the contact form $\alpha$ on $S$ can be extended to a 1-form
$\tilde{\alpha}$ in a neighborhood of $S\subset (X,\omega)$ in which $d\tilde{\alpha}=\omega$. Define the Liouville vector field $\Theta$ in the neighborhood by $i(\Theta)d\tilde{\alpha}=\tilde{\alpha}$. Then $\Theta$ is nowhere zero on $S$ and transversal to $S$ since $\alpha$ is a contact form on $S$. Let $\psi_t$ be the flow generated by $\Theta$, then $\psi_t(S)$,  $t\in (-\delta,\delta)$, parameterizes a neighborhood $(-\delta,\delta)\times S$ of $S$ in $X$ for some $\delta>0$.
On the other hand, $L_{\Theta}\tilde{\alpha}=d(i(\Theta)\tilde{\alpha})+
i(\Theta)d\tilde{\alpha}=
\tilde{\alpha}$, from which it follows that $\psi_t^\ast(\tilde{\alpha})=e^t\tilde{\alpha}$. So $\omega
=d\tilde{\alpha}=d(e^t\alpha)$, observing that $\tilde{\alpha}=\alpha$ on $S$
since $\Theta$ is transversal to $S$ and $\tilde{\alpha}(\Theta)=d\tilde{\alpha}(\Theta,\Theta)=0$.
\hfill $\Box$

\vspace{3mm}

\noindent{\bf The Proof of Theorem 1.2:}

The proof consists of three steps.

{\bf Step 1:} Let $S\subset (X,\omega)$ be a compact hypersurface of contact type and $\alpha$ be an induced contact form on $S$. By Lemma 2.3, $\omega=d(e^t\alpha)$ in a collar neighborhood $(-\delta,\delta)\times S$ of $S$ in $(X,\omega)$.
Let $X_1$ be the 4-manifold obtained by cutting $X$ open along $S$ and then inserting the cylinder $[0,1]\times S$ into it. Then $X_1$ is also a symplectic manifold with the symplectic form $\omega_1$ given by
$$ 
\omega_1=\left\{\begin{array}{cc}
\omega & \hspace{2mm} on \hspace{2mm} X\setminus (-\delta,\delta)\times S\\
d(\rho(t)\alpha) & \hspace{2mm} on \hspace{2mm} (-\delta, 1+\delta)\times S.
\end{array} \right.
$$
Here $\rho(t)$ is a smooth function on $(-\delta, 1+\delta)$ such that
$\rho^\prime(t)>0$ and $\rho(t)=e^t$ on $(-\delta,-\frac{1}{2}\delta)$ and
$\rho(t)=e^{t-1}$ on $(1+\frac{1}{2}\delta,1+\delta)$, and 
$\rho(0)=a$, $\rho(1)=b$ with $\rho(t)$ being linear on $[0,1]$ of slope $\epsilon_0$ where $a,b$ are in
$(e^{-\frac{1}{2}\delta}, e^{\frac{1}{2}\delta})$ with $ b-a=\epsilon_0>0$ (see
Figure 1, (a)).
Let $\sigma_l: (-\delta, l+\delta)\rightarrow (-\delta, 1+\delta)$ be a strictly increasing smooth function which equals to $t$ on $(-\delta,-\frac{1}{2}\delta)$ and $t-l+1$ on $(l+\frac{1}{2}\delta,l+\delta)$,
and maps $[0,l]$ linearly onto $[0,1]$ (see Figure 1, (b)). Let $X_l$ be the manifold obtained by
cutting $X$ open along $S$ and inserting the cylinder $[0,l]\times S$ into it.
Define a diffeomorphism $g_l: X_l\rightarrow X_1$ such that $g_l=(\sigma_l,id)$ on the neck and $g_l$ equals to identity on the rest. Then $X_l$ is a symplectic manifold with the pull back symplectic form $\omega_l=g^\ast_l\omega_1$. Note that the canonical class $c_1(K_l)$ of $(X_l,\omega_l)$ is the pull back of the canonical class $c_1(K_1)$ of $(X_1,\omega_1)$ via $g_l$.

\begin{figure}[h]
\centerline{\psfig{figure=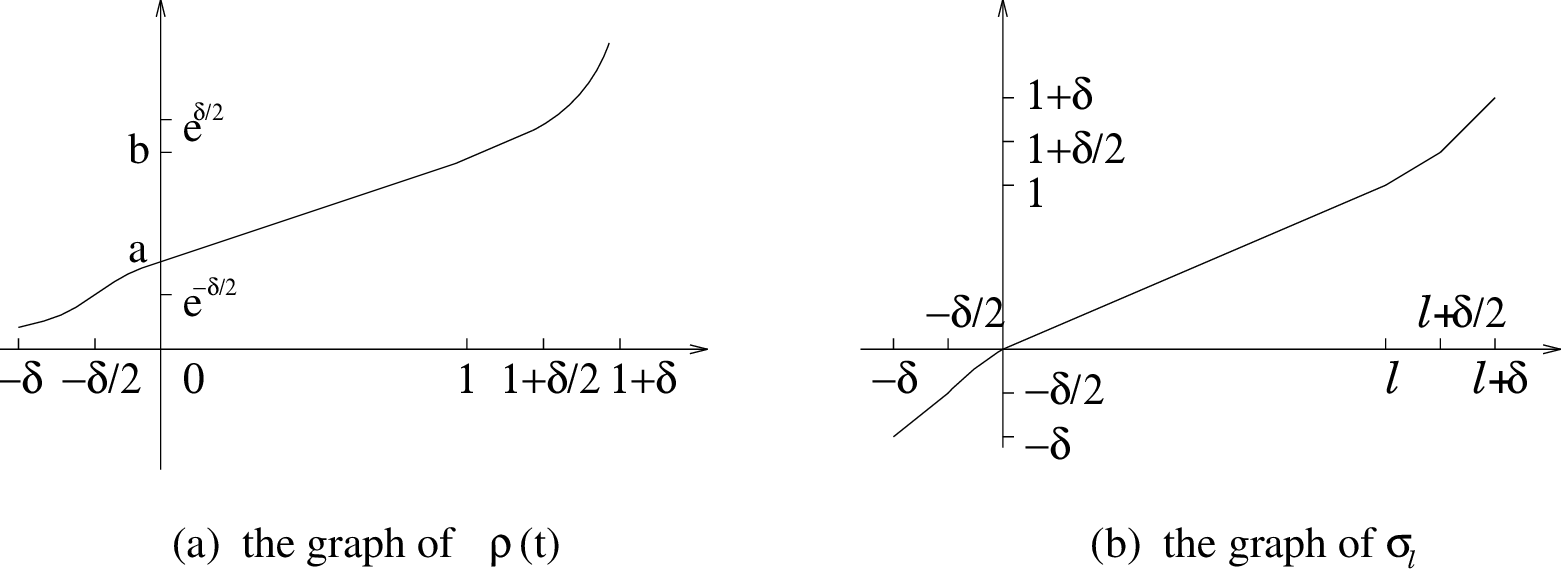,height=4.5cm,width=5in}}
\centerline{\large{Figure 1}}
\end{figure}

Recall that the Reeb vector field $v$ on $S$ associated to the contact form $\alpha$ is defined by $i(v)d\alpha=0$ and $\alpha(v)=1$. We define an almost complex structure $J_0$ on $\R\times S$ as follows:
$J_0$ equals to a fixed almost complex structure on $\xi=\{\alpha=0\}$ which is compatible to the symplectic form $d\alpha|_{\xi}$, and  $J_0(\frac{\partial}{\partial t})=v$, $J_0(v)=-\frac{\partial}{\partial t}$.

We will use a slightly extended version of Theorem 1.4 which was communicated
to the author by Taubes (see [McS]). In fact, in Taubes' theorem, the generic
compatible almost complex structures can be chosen freely in a region as 
long as no pseudo-holomorphic curve is supported in that region. On each symplectic manifold $(X_l,\omega_l)$, the neck $(-\delta,\delta+l)\times S$  does not support any pseudo-holomorphic curve since the symplectic form $\omega_l$ is exact on it. Therefore by the extended version of Taubes' theorem, for each $l$ there is a generic $\omega_l$-compatible almost complex structure $J_l$ on $(X_l,\omega_l)$ which equals to $J_0$ on the neck $[0,l]\times S$, and as $l\rightarrow\infty$, $\{J_l\}$ converges over $X_l\setminus (-\delta,\delta+l)\times S$ to a fixed $\omega$-compatible almost complex structure $J$ on $X$, and there is an embedded $J_l$-holomorphic curve $\Sigma_l$ in $(X_l,\omega_l)$
(may have several components) such that $c_1(K_l)$ is Poincare dual to $[\Sigma_l]$. 

\begin{lem}
The genus of each component of $\Sigma_l$ is uniformly bounded from above by
$c_1(K_1)^2+b^-_2(X_1)+1$.
\end{lem}

\noindent{\bf Proof:}
Suppose that $\{\Sigma^i_l\}$ are the components of $\Sigma_l$ and $[\Sigma^i_l]$ is Poincare dual to $e_i\in H^2(X_l,\Z)$. Then $c_1(K_l)=\sum e_i$. By the adjunction equality,
$$
2g_i-2=e_i^2+c_1(K_l)\cdot e_i=2e_i^2,
$$
where $g_i$ is the genus of $\Sigma_l^i$. One easily sees that if $e^2_i<0$,
then $e^2_i=-1$ and $g_i=0$. Therefore
$$
\sum e^2_i= c_1(K_l)^2+s_l\leq c_1(K_l)^2+b^-_2(X_l)=c_1(K_1)^2+b^-_2(X_1)
$$
where the summation is taken over all the $e_i's$ with $e_i^2\geq 0$, and $s_l$ is the number of $-1$-sphere components of $\Sigma_l$. Now it is easy to see that $g_i=1+e_i^2$ is bounded from above by $c_1(K_1)^2+b^-_2(X_1)+1$, which is
uniform in $l$.
\hfill $\Box$

\vspace{2mm}

{\bf Step 2:} 
Take a sequence of $l_n\rightarrow\infty$. For simplicity, the subscript $l_n$ is replaced by $n$ in the notation. Let $h_n$ be the hermitian metric on
$(X_n,\omega_n, J_n)$ defined by $h_n=\omega_n(\cdot,J_n\cdot)$, and $h_0$ be the hermitian metric on $\R\times S$ defined by $h_0=\omega_0(\cdot,J_0\cdot)$ where $\omega_0=dt\wedge \alpha+d\alpha$. Observe that the area of $\Sigma_n$ in $(X_n,J_n,h_n)$ is uniformly bounded from above by a constant $c( =\int_{\Sigma_1}\omega_1)$. Therefore, on each neck
$[0,l_n]\times S$, there is a sub-neck $I_n \times S$
(we can assume that the endpoints of $I_n$ are regular values of the function
$t_n^\prime=t|_{\Sigma_n}$) with the following property: $|I_n|=\delta_0$ for a small $\delta_0>0$ and $Area_{h_n}(\Sigma_n\bigcap (I_n\times S))$ is bounded from above by $\frac{1}{2}l_n^{-1}\epsilon_0\epsilon$ (recall that $\epsilon_0$ is involved in the definition of $\rho(t)$ in {\bf Step 1}). Here $\epsilon$ is chosen as in Lemma ~2.1 for the hermitian metric $h_0$ on $(\R\times S,J_0)$ and the annulus $[0,1]\times S^1$ ($S^1$ has unit length) with the standard complex structure and metric. Note that $Area_{h_0}(\Sigma_n\bigcap (I_n\times S))$ is bounded from above by $\frac{1}{2}\epsilon$ for large enough $n$. This is because on each neck $[0,l_n]\times S$,  
$h_0\leq (l_n\epsilon_0^{-1})h_n$ for large enough $n$. For each $n$, choose a sub-interval $I_n^1$ of $I_n$ such that 
the endpoints of $I_n^1$ are regular values of $t_n^\prime$, and $|I_n^1|=\frac{1}{2}|I_n|=\frac{1}{2}\delta_0$ and $dist(\partial I_n,\partial I_n^1)\geq\frac{1}{8}\delta_0$. 

\vspace{3mm}

We observe that the contact structure $\xi=\{\alpha=0\}$ as a complex line
bundle on $S$
is isomorphic to $K^{-1}|_S$ where $K$ is the canonical bundle of $(X,\omega)$. So for any regular value $t_0\in I_n^1$ of the function $t_n^\prime=t|_{\Sigma_n}$, $c_1(\xi)$ is Poincare dual to $\Sigma_n\bigcap(\{t_0\}\times S)$. Each  component of $\Sigma_n\bigcap(\{t_0\}\times S)$ is oriented by $-\alpha$ if we orient $S$ by $\alpha\wedge d\alpha$ (note that $t_0$ is a regular value
of $t_n^\prime$ so that the pull back of $\alpha$ to $\Sigma_n\bigcap(\{t_0\}\times S)$ is nowhere zero).

\begin{lem}
For any regular value $t_0\in I_n^1$, each component of $\Sigma_n\bigcap(\{t_0\}\times S)$ is connected through
$\Sigma_n\bigcap([t_0,l_n+\delta)\times S)$ to the inside of $X_n\setminus (-\delta, l_n+\delta)\times S$ (see Figure 2).
\end{lem}

\noindent{\bf Proof:}
If not, then it is easy to see that there are some components $\{\gamma_i\}$ of 
$\Sigma_n\bigcap(\{t_0\}\times S)$ which bound a pseudo-holomorphic curve $F\subset \Sigma_n\bigcap([t_0,l_n+\delta)\times S)$. If we orient $\{\gamma_i\}$ canonically as the boundary of the oriented surface $F$, then each integral $\int_{\gamma_i}\rho_n(t_0)\alpha$ is negative, where  $\rho_n(t)=\rho\circ\sigma_n(t)$ (see {\bf Step 1}). The reason is that $t_0$ is a regular value of $t_n^\prime=t|_{\Sigma_n}$ and $J_n(\frac{\partial}{\partial t})=v$ so that the pull back of $\alpha$ to $\gamma_i$ is a negative 
multiple of the volume form on $\gamma_i$. On the other hand,
$0<\int_{F}\omega_n=\sum_{i}\int_{\gamma_i}\rho_n(t_0)\alpha$ by the Stokes'
theorem. This contradiction proves the lemma.
\hfill $\Box$

\begin{figure}[h]
\centerline{\psfig{figure=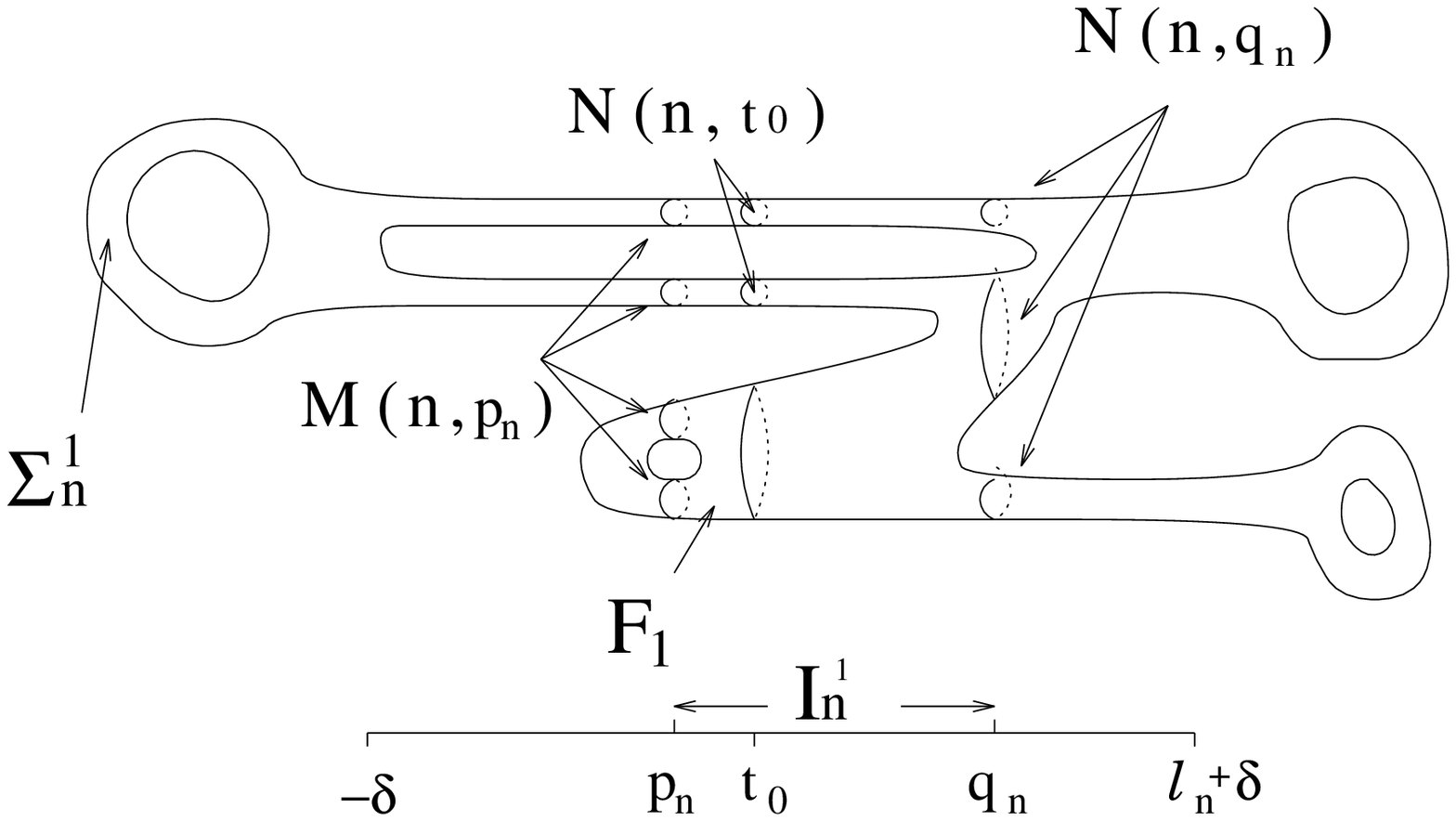,height=5cm,width=3.5in}}
\centerline{\large{Figure 2}}
\end{figure}

Now by Lemma 2.5, the components of $\Sigma_n\bigcap(\{t_0\}\times S)$ are divided into two groups:

\begin{enumerate}
\item A component of $\Sigma_n\bigcap(\{t_0\}\times S)$ belongs to Group {\bf I} if it is connected to the inside of $X_n\setminus (-\delta, l_n+\delta)\times S$ through both $\Sigma_n\bigcap((-\delta,t_0]\times S)$ and $\Sigma_n\bigcap([t_0,l_n+\delta)\times S)$. 
\item The rest of  $\Sigma_n\bigcap(\{t_0\}\times S)$ belongs to Group {\bf II} which bounds a pseudo-holomorphic 
curve $F_1\subset\Sigma_n\bigcap((-\delta,t_0]\times S)$ (may have several components). The homology class of the union of these components (each of which is oriented by $-\alpha$) is zero in $S$ (see Figure 2). 
\end{enumerate}

\begin{lem}
For any $n$ and regular value $t_0\in I_n^1$ of $t_n^\prime$,
there is a subset $N(n,t_0)$ of $\Sigma_n\bigcap(\{t_0\}\times S)$ such that $c_1(\xi)$ is Poincare dual to $N(n,t_0)$ and 
$\# N(n,t_0)\leq N$ for a constant $N$ independent of $n$ and $t_0\in I_n^1$.
\end{lem}

\noindent{\bf Proof:}
We define $N(n,t_0)$ to be the set of components of $\Sigma_n\bigcap(\{t_0\}\times S)$ which belong to Group {\bf I} (see Figure 2).
It remains to prove that there is an $N$ independent of $n$ and $t_0\in I_n^1$ such that $\# N(n,t_0)\leq N$. But this follows from the following two reasons:

\begin{itemize} 
\item the $h_n$-area of any component of $\Sigma_n\bigcap (X_n\setminus (-\delta, l_n+\delta)\times S)$ is bounded from below by a constant $c_1>0$
by Lemma 2.2 (Monotonicity), so the number of components of $\Sigma_n\bigcap (X_n\setminus (-\delta, l_n+\delta)\times S)$ is uniformly
bounded from above; 
\item the genus of each component of $\Sigma_n$ is uniformly bounded from above. \end{itemize}
\hfill $\Box$

Let $p_n$ and $q_n$ be the left and right endpoints of $I_n^1$. We define
a subset $M(n,p_n)$ of $\Sigma_n\bigcap(\{p_n\}\times S)$ as follows: a component $\gamma$ of $\Sigma_n\bigcap(\{p_n\}\times S)$ is in $M(n,p_n)$ if
there is a smooth path $\Gamma$ in $\Sigma_n\bigcap([-\delta, q_n]\times S)$ joining $N(n,q_n)$ to $\Sigma_n\bigcap(\{-\delta\}\times S)$ which intersects $\Sigma_n\bigcap(\{p_n\}\times S)$ first time at $\gamma$ geometrically once.
(See Figure 2.)

\begin{lem}
There is an $M$ independent of $n$ such that $\#M(n,p_n)\leq M$.
\end{lem}

\noindent{\bf Proof:}
It is easy to see from Lemma 2.5 that any element in $N(n,p_n)$ is connected to
$N(n,q_n)$ through $\Sigma_n\bigcap([p_n,q_n]\times S)$. So $N(n,p_n)\subset M(n,p_n)$. It suffices
to prove that $\#(M(n,p_n)\setminus N(n,p_n))$ has a uniform upper bound.

Let $\gamma$ be in $M(n,p_n)\setminus N(n,p_n)$. By the definition of $M(n,p_n)$, there is a smooth path $\Gamma$ in $\Sigma_n\bigcap([-\delta, q_n]\times S)$ joining $N(n,q_n)$ to $\Sigma_n\bigcap(\{-\delta\}\times S)$
which intersects $\Sigma_n\bigcap(\{p_n\}\times S)$ first time at $\gamma$ geometrically once. Since $\gamma$ is not in $N(n,p_n)$, $\Gamma$ must intersect
$\Sigma_n\bigcap(\{p_n\}\times S)$ at another component $\gamma^\prime$. By Lemma 2.5, $\gamma^\prime$ is connected to $X_n\setminus (-\delta,l_n+\delta)\times S$
through a path in $\Sigma_n\bigcap([p_n,l_n+\delta)\times S)$ which must 
intersect $\Sigma_n\bigcap(\{q_n\}\times S)$ at a component in $N(n,q_n)$. So
there is a path $\Gamma^\prime$ on $\Sigma_n$ with two ends in
$N(n,q_n)$ which intersects $\Sigma_n\bigcap(\{p_n\}\times S)$ only at $\gamma$
and $\gamma^\prime$. The relative homology class $[\Gamma^\prime]$ is non-zero in $H_1(\Sigma_n, N(n,q_n))$ since $\Gamma^\prime$ intersects 
$\gamma$ geometrically once. Now it follows from Lemmas 2.4 and 2.6 that $\#(M(n,p_n)\setminus N(n,p_n))$ has a uniform upper bound.
\hfill $\Box$

\vspace{3mm}

Now we pick a Morse function $t_n$ on $\Sigma_n\bigcap(I_n\times S)$ such that
$\|t_n-t_n^\prime\|_{C^k}<e^{-l_n}$ (recall that $t^\prime_n=t|_{\Sigma_n}$) 
and any two different critical points of $t_n$ have different values. We 
further require that $t_n=t^\prime_n$ at the endpoints of both $I_n$ and $I_n^1$. 

We associate a graph $G_n$ to each curve $\Sigma_n\bigcap (I_n\times S)$ by
the Morse function $t_n$ (see Figure 3). The graph $G_n$ has the following properties: 

\begin{itemize}
\item [{(a)}] each edge of $G_n$ corresponds to an open annulus in $\Sigma_n\bigcap (I_n\times S)$ on which $t_n$ is regular and each vertex corresponds to a critical point of $t_n$; 
\item [{(b)}] there is a projection $\pi_n:G_n\rightarrow I_n$ such that $\pi_n=t_n$ under the correspondence between $G_n$ and $\Sigma_n\bigcap (I_n\times S)$ ($\pi_n$ is one to one on each edge and maps the vertices to the critical values of $t_n$); 
\item [{(c)}] vertices correspondent to critical points of index $1$ have three edges (one or two on the left) and vertices correspondent to critical points of index $0$ (or $2$) have only one edge which is on the right (or left). 
\end{itemize}

Let $G_n^1$ be the sub-graph of $G_n$ which corresponds to $\Sigma_n\bigcap (I_n^1\times S)$, i.e. $G_n^1=\pi_n^{-1}(I_n^1)$.

\begin{figure}[h]
\centerline{\psfig{figure=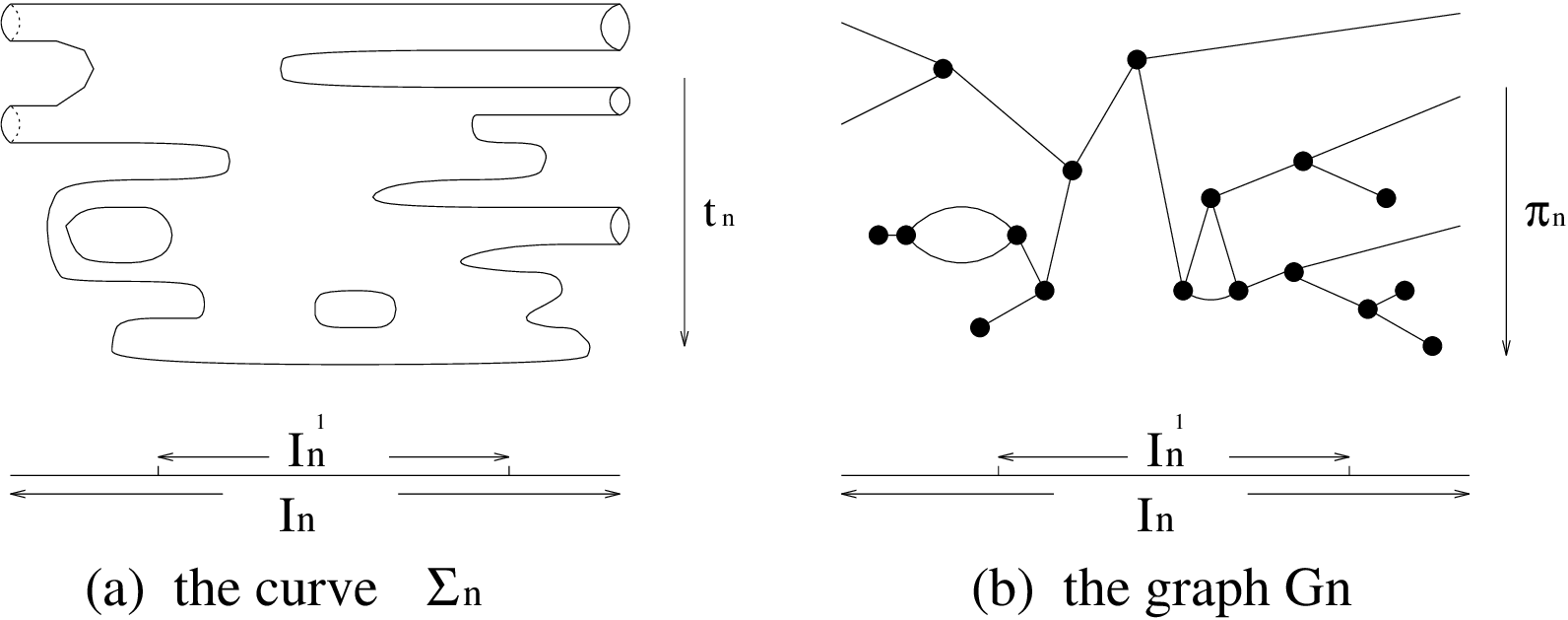,height=4.5cm,width=5in}}
\centerline{\large{Figure 3}}
\end{figure}

Let $\Lambda_n$ be the sub-graph of $G^1_n$ which is the union of all the paths in $G_n^1$ joining the points correspondent to $M(n,p_n)$ with the points 
correspondent to $N(n,q_n)$ on $G_n^1$. Clearly the number of paths in $\Lambda_n$ has
a uniform upper bound since $\#M(n,p_n)\leq M$, $\#N(n,q_n)\leq N$ and 
the genus of each component of $\Sigma_n$ has a uniform upper bound.

\begin{lem}
There exists a $\delta_1>0$ such that
for large $n$, $\Lambda_n$ has the following property:
\begin{itemize} 
\item there are annuli $A_{n,i}$ in $\Sigma_n\bigcap(I_n\times S)$ such that
$\partial A_{n,i}$ are regular level circles of $t_n^\prime$ with regular 
values $x_n, y_n$ ( the same for all i) satisfying $y_n-x_n=\delta_1$;
\item $N(n,x_n)$ consists of the left boundary components of $A_{n,i}$, in particular, $\#A_{n,i}\leq N$.
\end{itemize}
\end{lem}

\noindent{\bf Proof:}
Let $\Gamma$ be a path in $G_n^1$ such that the projection $\pi_n|_\Gamma:
\Gamma\rightarrow I_n^1$ is surjective.
We put the vertices on $\Gamma$ (critical points of $t_n$ of Morse index $1$) into four groups as follows:

\begin{enumerate}
\item Group $(1)$ consists of those vertices whose only edge that is not on the path $\Gamma$ is part of a closed cycle in $G_n$ involving part of $\Gamma$. Each
of such vertices is associated with a closed cycle in $G_n$ such that the vertex
is the most left or the most right one amongst the vertices on both $\Gamma$ and the closed cycle. 
It is easy to see that the number of vertices in Group $(1)$ is at most twice
of the number of the associated closed cycles.
\item Group $(2)$ consists of those vertices which start a sub-graph of $G_n$ (never
connect to $\Gamma$ through another vertex) such that there is a closed cycle in this sub-graph.
\item Group $(3)$ consists of those vertices which start a tree that
goes all the way to one of the two ends of $G_n$
(i.e. $\pi_n^{-1}(\partial I_n)$). 
\item Group $(4)$ consists of those vertices
which start a tree whose projection under $\pi_n$ is within $I_n$.
\end{enumerate}

 ( In Figure 4, (b), the path $\Gamma$ is thickened. The vertex $p_1$ is in 
Group $(2)$, $p_2$ in Group $(3)$, $p_3, p_4$ in Group $(1)$ and $p_5$ in Group $(4)$.)

Since the genus of each component 
of $\Sigma_n$ is uniformly bounded from above,  
the number of elements in Groups $(1),(2)$ has an upper bound independent of $n$ or the path $\Gamma$. The
number of elements in Group $(3)$ is also
uniformly bounded from above. The reason is that each component of the pseudo-holomorphic
curve in $(I_n\setminus int{I_n^1})\times S$ with boundary on both of the
boundaries of $I_n\times S$ and $I_n^1\times S$ has $h_0$-area bounded
from below by a constant $c_2>0$ (independent of $n$ or $\Gamma$) by Lemma 2.2 (Monotonicity). Therefore the number of vertices on $\Gamma$ which belong to Groups $(1)$, $(2)$ or $(3)$ is uniformly bounded from above. 

Since the number of paths in $\Lambda_n$ is uniformly bounded from above, the 
total number of vertices on $\Lambda_n$ which belong to Groups $(1)$, $(2)$ or $(3)$ for some path in $\Lambda_n$ is also uniformly bounded from above.
Let $M^\prime$ be such an upper bound, and set $\delta_1=(4M^\prime+8)^{-1}\delta_0$.
Then there is an interval $I_n^2$ in $I^1_n$ whose endpoints are regular values of $t_n$ such that $|I_n^2|=(M^\prime+2)^{-1}|I_n^1|=(2M^\prime+4)^{-1}\delta_0=2\delta_1$
and $I_n^2$ does not contain $\pi_n(p)$ for any vertex $p$ on $\Lambda_n$
which belongs to Groups $(1)$, $(2)$ or $(3)$ for some path in $\Lambda_n$ (but $I_n^2$ may contain $\pi_n(p)$ for some $p$ in Group $(4)$). Now look at 
$\pi_n^{-1}(\partial I_n^2)\bigcap\Lambda_n$. 
It is easy to see that the corresponding circles on $\Sigma_n\bigcap (I_n^1\times S)$ bound a set of annuli in $\Sigma_n\bigcap (I_n\times S)$, 
since the tree started by any vertex in Group $(4)$ (for example, the vertex $p_5$ in Figure 4, $(b)$) corresponds to a topological disc in $\Sigma_n\bigcap (I_n\times S)$. (The case when $\Lambda_n$ consists of only one path is shown 
in Figure 4).

\begin{figure}[h]
\centerline{\psfig{figure=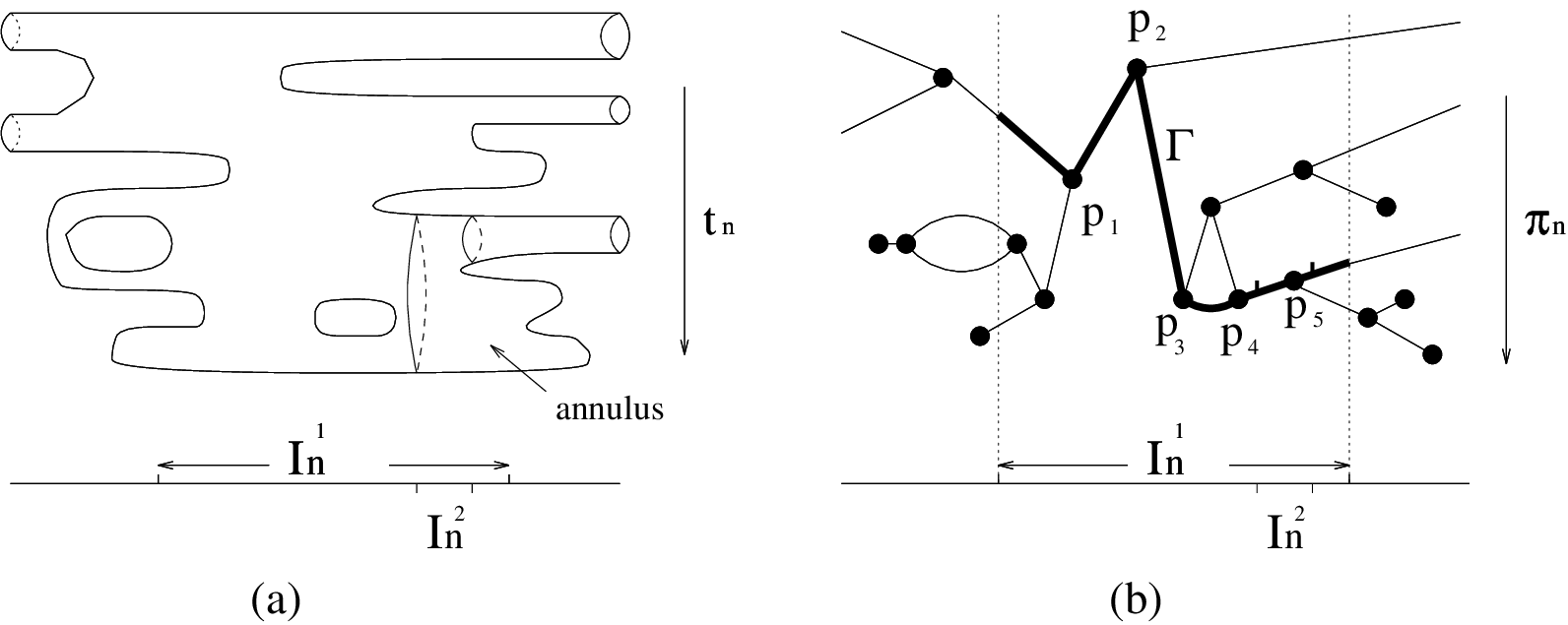,height=4.5cm,width=5in}}
\centerline{\large{Figure 4}}
\end{figure}

Pick regular values $x_n$, $y_n$ ($x_n<y_n$) of $t_n^\prime$ in $I_n^2$ such that $|x_n-y_n|=\delta_1(=\frac{1}{2}|I_n^2|)$ and $\min(dist(x_n,\partial I_n^2), dist(y_n,\partial I_n^2))>\frac{1}{8}\delta_1$.
By the construction of $\Lambda_n$,
$N(n,x_n)$ must lie in the part of $\Sigma_n$ that corresponds to $\Lambda_n$,
and for large $n$, must lie in the annuli constructed in the previous paragraph 
(since $\|t_n-t_n^\prime\|_{C^k}<e^{-l_n}$). Each annulus contains at most one
component of $N(n,x_n)$ which has non-zero homology in that annulus. 
Collect all the annuli that contain a component of $N(n,x_n)$. Each of these annuli also contains one component of $N(n,y_n)$ so that it has a sub-annulus
which is bounded by a component of $N(n,x_n)$ on the left and a component
of $N(n,y_n)$ on the right. These sub-annuli are the claimed annuli $A_{n,i}$. 
\hfill $\Box$

\begin{lem}
There exist sequences of pseudo-holomorphic
maps $f_{n,i}:[0,L]\times S^1\rightarrow I_n\times S$ for an $L\in (0,1]$ such that the image of each $f_{n,i}$ lies in the annulus $A_{n,i}$ (constructed
in Lemma 2.8) and $f_{n,i}(\{0\}\times S^1)$ is the left boundary component of $A_{n,i}$. After applying a translation to $f_{n,i}$ if necessary, $\{f_{n,i}\}$ has a subsequence which is $C^\infty$ convergent on
$[\frac{L}{4},\frac{3L}{4}]\times S^1$. Furthermore, the $S$-component of the image of each limiting pseudo-holomorphic map $f_i:[\frac{L}{4},\frac{3L}{4}]\times S^1\rightarrow \R\times S$ lies in a Reeb orbit.
\end{lem}

\noindent{\bf Proof:} 
First recall a fact from the classical theory of Riemann surfaces. Let $A$ be an annulus equipped with a smooth Riemannian metric $\mu$. The complex structure determined by $\mu$ has a uniformization. More precisely, $A$ is holomorphically equivalent to the product $[0,L]\times S^1$ with the standard complex structure and metric ($S^1$ has unit length), where $L$ is the modulus of $A$, and is given by
$$
{L^{-1}}=\inf\int_{A}|du|^2
$$
over smooth real valued functions $u$ on $A$ which take the value $0$ on one boundary component and $1$ on the other (see [P]). 

For our purpose here, we need to prove that the modulus $L_{n,i}$ of each
annulus $A_{n,i}$ has a lower bound independent of $n$ and $i$. This goes
as follows. Recall that $t_n^\prime=x_n$ on the
left boundary component of each $A_{n,i}$ and $t_n^\prime=y_n$ on the right boundary component. Define a smooth function $u_{n,i}$ on $A_{n,i}$ by $u_{n,i}=\delta_1^{-1}(t_n^\prime-x_n)$ ($\delta_1=y_n-x_n$). Then we have
$$
{L_{n,i}^{-1}}\leq \int_{A_{n,i}}|du_{n,i}|^2\leq \delta_1^{-2}\int_{A_{n,i}}|dt_n^\prime|^2
\leq 2\delta_1^{-2}Area_{h_0}(A_{n,i})\leq \delta_1^{-2}\epsilon.
$$
So $L_{n,i}\geq \epsilon^{-1}\delta_1^2$.  Now it is easy to see that we have  sequences of pseudo-holomorphic maps $f_{n,i}:[0,L]\times S^1\rightarrow I_n\times S$ with $L=\min(\epsilon^{-1}\delta_1^2,1)$ such that the image
of $f_{n,i}$ lies in $A_{n,i}$ and $f_{n,i}(\{0\}\times S^1)$ is the left boundary component of $A_{n,i}$. Equip $I_n\times S$ with the hermitian metric $h_0$. Then
the energy $E(f_{n,i})$ is bounded from above by $\frac{1}{2}\epsilon$.
By Lemma ~2.1, the $C^0$ norm of $df_{n,i}$ over $[\frac{L}{4},\frac{3L}{4}]\times S^1$ is uniformly bounded from above. By the regularity theory of pseudo-holomorphic maps, after applying a translation to $f_{n,i}$ if necessary, $\{f_{n,i}\}$ has a subsequence which is $C^\infty$ convergent on $[\frac{L}{4},\frac{3L}{4}]\times S^1$. Let $f_i:[\frac{L}{4},\frac{3L}{4}]\times S^1\rightarrow \R\times S$ denote the limit
(here $\R\times S$ is equipped with the almost complex structure $J_0$, 
see {\bf Step 1}). Then the $S$-component of the image of $f_i$ lies in a Reeb orbit. This is because 
$$
\int_{[\frac{L}{4},\frac{3L}{4}]\times S^1} f^\ast_{n,i}(d\alpha)\leq e^{\delta} Area_{h_n}(A_{n,i})
\leq e^{\delta}l_n^{-1}\epsilon_0\epsilon
$$ 
which goes to $0$ as $n\rightarrow\infty$.
\hfill $\Box$

\vspace{3mm}

Recall that $N(n,x_n)=\{f_{n,i}(\{0\}\times S)\}$, which is Poincare dual to $c_1(\xi)$, is oriented by $-\alpha$, where $\alpha$ is the contact form on $S$ and $\xi=\{\alpha=0\}$ is the contact structure.  We orient $\{\frac{L}{2}\}\times S^1$ as the negative of the oriented boundary component of oriented
surface $[0,\frac{L}{2}]\times S^1$ so that 
$\bigcup_{i} f_{n,i}(\{\frac{L}{2}\}\times S^1)$
is homologous to $N(n,x_n)$. Let $\tilde{f_i}$ be the projection of $f_i$ into $S$. Then $c_1(\xi)$ is Poincare dual to $\bigcup_i \tilde{f_i}(\{\frac{L}{2}\}\times S^1)$ with each $\tilde{f_i}(\{\frac{L}{2}\}\times S^1)$ lying in a Reeb orbit. We can throw
away any $\tilde{f_i}(\{\frac{L}{2}\}\times S^1)$ which is not a closed orbit, since its homology class is zero. So to finish the 
proof of the first assertion in Theorem 1.2, it suffices to show that for each
$i$, $\int_{\tilde{f_i}(\{\frac{L}{2}\}\times S^1)}\alpha\leq 0$. But this follows from 
$$
\int_{f_{n,i}(\{\frac{L}{2}\}\times S^1)}(\rho\circ\sigma_n)\alpha
=\int_{f_{n,i}(\{0\}\times S^1)}(\rho\circ\sigma_n)\alpha -
\int_{f_{n,i}([0,\frac{L}{2}]\times S^1)}d((\rho\circ\sigma_n)\alpha)<0.
$$
(Recall from {\bf Step 1} that $\omega_n=d((\rho\circ\sigma_n)\alpha)$ on the neck $[0,l_n]\times S$.)

\vspace{3mm}

{\bf Step 3:}
Suppose that $S$ bounds a submanifold $W\subset (X,\omega)$ such that $c_1(W)\neq 0$ and $\omega$ is exact on $W$. To finish the proof we need to 
show that $S$ has at least one closed characteristic. 

Let $\omega=d\lambda$ on $W$. Then one can extend $\lambda$ to a neighborhood of $W\subset X$ in which $\omega=d\lambda$ still holds ([W]). We can assume that
$(-\delta,\delta)\times S$ is contained in this neighborhood. For each $l_n\geq 1$, let $W_n=W\bigcup([0,l_n+\delta)\times S)$ in $(X_n,\omega_n)$. Then there is a 1-form $\lambda_n$ on $W_n$ such
that $\omega_n=d\lambda_n$. The proof is as follows: first let $\lambda_1=g^\ast\lambda$ on $W_1$ where $g: X_1\rightarrow X$ is a diffeomorphism which is given by
$(t,x)\rightarrow (\ln\rho(t),x)$ on the neck, then let $\lambda_n=g^\ast_n\lambda_1$ where the diffeomorphism $g_n: X_n\rightarrow X_1$ is defined in {\bf Step ~1}.
Without loss of generality, we assume that the
normal vector field $\frac{\partial}{\partial t}$ is outward with respect to $W$.

Note that there is at least one component of each $\Sigma_n$ which intersects  both of $W\subset X_n$ and $X_n\setminus W_n$ since $c_1(W)\neq 0$ and $\omega_n$ is exact on $W_n\subset X_n$. In particular,
$N(n,t_0)$ is not empty for any $n$ and regular value $t_0\in I_n^1$ of $t^\prime_n=t|_{\Sigma_n}$.

We recall from {\bf Step 2} that $t_n^\prime=x_n$ on the left boundary component of each annulus $A_{n,i}$ (for all $i$). By the definition of $N(n,x_n)$, the union of the left boundary components of annuli $A_{n,i}$ (which is 
$N(n,x_n)$) bounds a
pseudo-holomorphic curve $\Sigma_n^1$ in $W_n$ which goes through the inside of $W$ (see Figure 2). By Lemma 2.2, the $h_n$-area of $\Sigma_n^1$ is bounded from below by a constant $c_3>0$. So by Stokes'
theorem and the fact that $\omega_n=d\lambda_n$ in $W_n$, we have $\max_{\gamma\in\{\gamma_{n,i}\}}\int_{\gamma}\lambda_n\geq N^{-1}c_3>0$, where $\{\gamma_{n,i}\}$ are the left boundary components of $A_{n,i}$, each of which 
is oriented by $\alpha$.

Pick an annulus $A_n\in \{A_{n,i}\}$ whose left boundary component $\gamma_n\in \{\gamma_{n,i}\}$ satisfies 
$$
\int_{\gamma_{n}}\lambda_n\geq N^{-1}c_3>0.
$$
Let $f_n:[0,L]\times S^1\rightarrow I_n\times S$ be the sequence of pseudo-holomorphic maps associated to $A_n$, and $f:[\frac{L}{4},\frac{3L}{4}]\times S^1\rightarrow \R\times S$ be the limit of $f_n$ (see Lemma 2.9).
Then the $S$-component of the image of $f$ must lie in a closed Reeb orbit. The proof is as follows. Suppose it does not lie in a closed orbit. Let $\gamma=\{\frac{L}{2}\}\times S^1$ be the circle in $[0,L]\times S^1$ canonically oriented as a boundary component of $[0,\frac{L}{2}]\times S^1$, 
and $s$ be the time coordinate of the Reeb flow, then $\int_{\gamma} \tilde{f}^\ast \alpha=\int_{\gamma} \tilde{f}^\ast ds=\int_{\gamma} d\tilde{f}=0$ where $\tilde{f}$ is the projection of $f$ into  $S$. On the other hand, $\int_{\gamma}f_n^\ast\lambda_n=
\int_{\gamma}f^\ast_n((\rho\circ\sigma_n)\alpha)$ for large $n$, since the homology class of  $\tilde{f}(\gamma)$ is zero in $S$ and the difference $\lambda_n-(\rho\circ\sigma_n)\alpha$ is closed on the neck $[0,l_n]\times S$.
This leads to a contradiction as follows: on the one hand, $\int_{\gamma}f_n^\ast \lambda_n>\int_{\gamma_{n}}\lambda_n\geq N^{-1}c_3>0$, 
on the other hand, let
$T_n(1)=\min (t|_{f_n(\gamma)})$ and $T_n(2)=\max (t|_{f_n(\gamma)})$, and 
$\tilde{f_n}$ be the $S$-component of $f_n$, then
$$
|\int_{\gamma}f^\ast_n((\rho\circ\sigma_n)\alpha)|
\leq
\int_{\gamma}|\max(\rho\circ\sigma_n)^\prime|
|T_n(2)-T_n(1)||\tilde{f_n}^\ast\alpha|+ |\int_{\gamma}\rho\circ\sigma_n(T_n(1))\tilde{f_n}^\ast\alpha|
$$
which goes to zero as $n \rightarrow\infty$, since $\int_{\gamma}\tilde{f_n}^\ast\alpha$ and 
$(\rho\circ\sigma_n)^\prime=\epsilon_0 l_n^{-1}$
go to zero and $T_n(2)-T_n(1)$, $|\tilde{f_n}^\ast\alpha|$ and $\rho\circ\sigma_n(T_n(1))$ remain bounded
as $n\rightarrow \infty$.
Therefore Theorem 1.2 is proved.
\hfill $\Box$

\vspace{2mm}

{\Small Current Address: UW-Madison, Madison, WI 53706, USA, 
{\it e-mail:} wechen@@math.wisc.edu}

\end{document}

%% file: def.tex
\let\text=\mbox
\newcommand{\skipthistext}[1]{}

\newcommand{\R}{{\mathbb R}}
\newcommand{\Z}{{\mathbb Z}}

\newcommand{\E}{\text{$\cal L$}}

\let\phi=\varphi
\let\epsilon=\varepsilon

\newcounter{thm}

\def\thmnumber{
\addtocounter{thm}{1}\if\arabic{section}0\relax
\else\arabic{section}.\fi%
\if\arabic{subsection}0\relax
\else\arabic{subsection}.\fi%
\arabic{thm}
\hspace{2mm}}

\def\setthmnumber{\setcounter{thm}{0}}

\newenvironment{thm}{
\vspace{3mm}\par\noindent 
{\bf Theorem \thmnumber} \it}
{\vspace{3mm}\par\rm}

\newenvironment{lem}{
\vspace{3mm}\par\noindent 
{\bf Lemma \thmnumber} \it}
{\vspace{3mm}\par\rm}

\newenvironment{coro}{
\vspace{3mm}\par\noindent 
{\bf Corollary \thmnumber}\it}
{\vspace{3mm}\par\rm}

\newenvironment{conj}{
\vspace{3mm}\par\noindent 
{\bf Conjecture \thmnumber}\it}
{\vspace{3mm}\par\rm}

\newcommand{\sectioni}[1]{\setthmnumber\setcounter{subsection}{0}\section{#1}}

\renewcommand{\subsection}[1]{%
\vspace{7mm}\par\noindent%
\addtocounter{subsection}{1}%
{\sc\arabic{section}.\arabic{subsection} \ {#1}}
\vspace{4mm}\par\noindent\setthmnumber}
